\documentclass[square]{ws-procs975x65}
\begin{document}
\title{Spectral properties of discrete alloy-type models}

\author{Martin Tautenhahn and Ivan Veseli\'c}

\address{Emmy-Noether-Projekt \emph{Schr\"o{}dingeroperatoren}, Fakult\"a{}t f\"u{}r Mathematik, \\ Technische Universit\"a{}t Chemnitz, 09127 Chemnitz, Germany \\
www.tu-chemnitz.de/mathematik/enp/}

\begin{abstract}
We discuss recent results on spectral properties of 
discrete alloy-type random Schr\"odinger operators.
They concern Wegner estimates and bounds on the fractional moments of the
Green's function.
\end{abstract}


\bodymatter
%
%
\section{Introduction and model} \label{sec:model}
A discrete alloy-type model is a family of operators $H_\omega=H_0+V_\omega$
on $\ell^2(\mathbb{Z}^d)$ where $H_0$ is the negative discrete Laplacian given by $(H_0
\psi)(x) = -\sum_{\lvert e \rvert = 1} \psi (x+e)$, and $V_\omega$ the
multiplication operator by the function
\begin{equation} \label{eq:alloy}
V_\omega (x) =\sum_{k \in \mathbb{Z}^d} \omega_k \, u(x-k) .
\end{equation}
Here $\omega_k, k \in \mathbb{Z}^d$, is an i.\,i.\,d.{} sequence 
of random variables and $u \colon \mathbb{Z}^d\to \mathbb{R}$ a so-called single-site
potential. The minimal conditions which are assumed throughout this note are
that the distribution of $\omega_0$ has a density $\rho \in L^\infty (\mathbb{R})$ and
that $u$ belongs to $\ell^1(\mathbb{Z}^d)$.
There are several well-studied relatives of the discrete alloy-type model:
\begin{romanlist}[(iii)]
\item The continuum alloy-type model $-\Delta +V_\omega$ on $L^2(\mathbb{R}^d)$ 
where $V_\omega $ has the same form as in \eqref{eq:alloy}. 
\item The Anderson model on $\ell^2(\mathbb{Z}^d)$: 
this is the special case of the discrete 
alloy-type model where $u=\delta_0$. 
\item The correlated Anderson model, which corresponds 
to $u=\delta_0$ and the case where the random coupling constants
$\omega_k$ are no longer assumed independent. %
\end{romanlist}
Obviously, the discrete alloy-type model is a special case of case (iii).
However, previous results on (iii)-type models require certain regularity
conditions 
which are typically not satisfied for alloy-type potentials, cf. Section \ref{sec:example}. 
\par
The main challenge of models of the type \eqref{eq:alloy}
is that $u$ may change its sign which leads to negative correlations 
of the potential. This results in a non-monotone dependence 
on the random variables of certain spectral quantities.
For this reason many established tools for the spectral 
analysis of random Schr\"odinger operators are not directly applicable
in this model. 
\par
For literature concerning the discussion in this section we refer to the references of \cite{elgart2009,veselic2009}.	

%
%
\section{Results}
Our results concern averages of resolvents and projectors. The presented results
are closely related to dynamical and spectral localisation. 
We use the symbol $\mathbb{E}$ to denote the average over the collection of random
variables $\omega_k$, $k \in \mathbb{Z}^d$. For $\Lambda \subset \mathbb{Z}^d$ we
denote by $H_{\omega , \Lambda} : \ell^2 (\Lambda) \to \ell^2 (\Lambda)$ the
natural restriction of $H_\omega$ to the set $\Lambda$.
%
%
%
\subsection{Wegner estimate}
We present here a selection of the results proven in \cite{veselic2009}.
For $L\in \mathbb{N}$ we denote the set $[0,L]^d \cap \mathbb{Z}^d$ by $\Lambda_L$. 
The number of points in the support of $u$ is denoted by $ \mathrm{rank}\,  u$.

\begin{theorem} \label{thm:degenerate}
Assume that the single-site potential $u$ has support in $\Lambda_{n}$ and
$\rho$ is of bounded variation. 
Then there exists a constant $c_u$ depending only on $u$ such that for any 
$L\in \mathbb{N}$, $E \in \mathbb{R}$ and $\epsilon>0$ we have
\[
\mathbb{E} \left \{ \mathrm{Tr}\,  \big
[\chi_{[E-\epsilon,E-\epsilon]}(H_{\omega,\Lambda_L})\big]\right\}
\le 
c_u \,   \|\rho\|_{\mathrm{Var}} \,   \mathrm{rank}\,  u  \ \epsilon\, (L+n)^{d\cdot (n+1)} .
\]
\end{theorem}
Here $\|\rho\|_{\mathrm{Var}} $ denotes the total variation of $\rho$.
By the assumption on the support of the single-site potential, 
$ \mathrm{rank}\,  u \le n^d$. Our bound is linear in the energy-interval length and
polynomial in the volume of the cube $\Lambda_L$. This implies that the Wegner bound can be
used for a localisation proof via multiscale analysis, as soon as an appropriate
initial length scale estimate is at disposal.
\begin{theorem} \label{thm:non-degenerate}
Assume  $\bar u := \sum_{k\in\mathbb{Z}^d} u(k) \neq 0$ and that $ \rho$ has compact
support. 
Let $m \in \mathbb{N} $ be such that $\sum_{\|k\|\ge m} |u(k)| \le |\bar u/2|$. 
Then we have
for any $L\in \mathbb{N}$, $E \in \mathbb{R}$ and $\epsilon>0$ 
\[
\mathbb{E} \left \{ \mathrm{Tr}\,  \big
[\chi_{[E-\epsilon,E-\epsilon]}(H_{\omega,\Lambda_L})\big]\right\}
\le 
\frac{8}{\bar u}  \, \|\rho\|_{\mathrm{Var}} \,  \min\big(L^d, \mathrm{rank}\,  u\big) \ \epsilon\,
(L+m)^{d} .
\]
\end{theorem}
In the case that the support of $u$ is compact we obtain a bound which is linear
in the 
volume of the box and thus yields the Lipschitz continuity of the integrated
density of states.
%
%
%
\subsection{Boundedness of fractional moments of Green's function}
In this subsection we assume that the function $\hat u \colon [0,2\pi)^d \to
\mathbb{C}$, defined by
\[
 \hat u (\theta) = \sum_{k \in \mathbb{Z}^d} u(k) {\rm e}^{\ensuremath{{\mathrm{i}}}k \cdot \theta} ,
\]
does not vanish, $ \mathrm{supp}\,  u$ is compact, and $\rho \in W^{1,1} (\mathbb{R})$. Let
$\Lambda \subset \mathbb{Z}^d$ be finite. For $x,y \in \Lambda$ and $z \in \mathbb{C} \setminus
\mathbb{R}$ we set $G_{\omega , \Lambda} (z;x,y) = \langle \delta_x ,
(H_{\omega, \Lambda} - z)^{-1} \delta_y \rangle$. 
%
\begin{theorem} \label{thm:apriori-trans}
  Let $s \in (0,1)$ and $\Lambda \subset \mathbb{Z}^d$ be finite. Then there exists a
constant
$C_u$ depending only on $u$, such that for all $z \in \mathbb{C} \setminus \mathbb{R}$ and all
$x,y \in \Lambda$
\[
 \mathbb{E} \Bigl\{ \bigl| G_{\omega , \Lambda} (z;x,y) \bigr|^s \Bigr\} \leq \bigl(C_u
\Vert\rho'\Vert_{L^1} \bigr)^s \frac{2^{1+s} s^{-s}}{1-s} .
\]
\end{theorem}
Theorem~\ref{thm:apriori-trans} is proven in Section \ref{sec:proof}. The proof
is a combination of ideas from \cite{tautenhahn2007} and \cite{veselic2008}. The
proof also gives rise to a quantitative estimate on the constant $C_u$. 
To our knowledge, Theorem \ref{thm:apriori-trans} does not immediately imply exponential decay of
fractional moments of the Green's function by standard methods.
%
%
%
\subsection{Exponential decay of fractional moments of Green's function}
In this subsection we assume $d = 1$ and $ \mathrm{supp}\,  u$ compact. 
Denote the diameter of the support by $n - 1 \in \mathbb{N}_0$. We define the
following conditions which may or may not hold:
\begin{itemize}
 \item[(A)] $ \mathrm{supp}\,  u = \{0,1,\dots,n-1\}$,
 \item[(B)] $ \mathrm{supp}\,  \rho$ compact,
 \item[(C)] $\rho \in W^{1,1} (\mathbb{R})$.
\end{itemize}
For $x,y \in \mathbb{Z}^d$ and $z \in \mathbb{C} \setminus \mathbb{R}$ we denote the Green's function
by $G_\omega (z;x,y) = \langle \delta_x , (H_\omega - z)^{-1} \delta_y \rangle$ 
and set 
\[
 C_{u,\rho} = \Bigl | \prod_{k \in \Theta} u(k) \Bigr
|^{-s/n}\Vert\rho\Vert_{\infty}^{s} \frac{2^{s} s^{-s}}{1-s} .
\]
\begin{theorem} \label{thm:expd1}
 Let $s \in (0,1)$ and assume that (A) holds. Then there exists a constant $C$,
depending on $\lVert \rho \rVert_\infty$, $s$ and $u$, such that for all $x,y
\in \mathbb{Z}$ with $|x-y| \geq n$ and all $z \in \mathbb{C} \setminus \mathbb{R}$ we have
\begin{equation} \label{eq:exp}
 \mathbb{E} \Bigl\{ \bigl| G_{\omega} (z;x,y) \bigr|^{s/n} \Bigr\} \leq C \,\, \exp \left\{-m \left \lfloor \frac{|x-y|}{n} \right\rfloor \right\} ,
\end{equation}
where $m = - \ln C_{u,\rho} $ and $\lfloor \cdot \rfloor$ is defined by $\lfloor t \rfloor = \max\{k\in \mathbb{Z} \colon k\leq t\}$.
\end{theorem}
If $C_{u,\rho} < 1$, or equivalently $m>0$, Ineq. \eqref{eq:exp} describes
exponential decay.
\begin{theorem} \label{thm:expd1_gen}
 Let $s \in (0,1/2)$ and assume that either (B) holds with
$\lVert\rho\rVert_\infty$ sufficiently small, \textbf{or} (C) holds with
$\lVert\rho'\rVert_{L^1}$ sufficiently small. Then there exist constants $m,C
\in (0,\infty)$, such that Eq. \eqref{eq:exp} holds for all $x,y \in \mathbb{Z}$ with
$|x-y| \geq 4n$ and all $z \in \mathbb{C} \setminus \mathbb{R}$.
\end{theorem}
Theorem~\ref{thm:expd1} and \ref{thm:expd1_gen} are proven in a joint paper \cite{elgart2009} with A.~Elgart.
The estimates of Theorem \ref{thm:expd1} and \ref{thm:expd1_gen} concern only
off-diagonal elements. If we assume (B) and $s\in (0,1/4n)$, then $\mathbb{E} \bigl\{|G_\omega (z;x,y)|^s
\bigr\}$ is uniformly bounded for $x,y \in \mathbb{Z}$ and $z \in \mathbb{C} \setminus \mathbb{R}$, see \cite[Appendix]{elgart2009}. In spite these estimates on the Green's function, neither dynamical
nor spectral localisation follow immediately using the existent methods.
\par
Currently, the analogues of Theorem~\ref{thm:expd1} and \ref{thm:expd1_gen} in higher dimension are an open question. In a forthcoming paper with A.~Elgart we intend to extend Theorem~\ref{thm:expd1_gen} to operators on a one-dimensional strip, i.\,e. $\Lambda = \Gamma \times \mathbb{Z}$ with $\Gamma \subset \mathbb{Z}^{d}$ finite.
%
%
%
\section{Example} \label{sec:example}
We show that discrete alloy-type models do not satisfy in general the regularity conditions required in \cite{aizenman1993,aizenman2001} for correlated Anderson models.
We consider $d = 1$. Let the density function $\rho$ be of bounded support
and let $u(0)=1$, $u(-1) = a \not = 0$ and $u(k) = 0$ for $k \in \mathbb{Z} \setminus \{-1,0\}$. For simplicity we require that the infimum of the support of $\rho$ is zero. Let $B_\epsilon$ and $A_\epsilon$ be the events $B_\epsilon = \{\omega \in \Omega : V_\omega(-1),V_\omega(1) \in [0,\epsilon]\}$ and $A_\epsilon = \{\omega \in \Omega : V_\omega(0) \in [0,\epsilon (a + \frac{1}{a})]\}$ with $\epsilon > 0$. Then, if $a>0$ one calculates $B_\epsilon \subset A_\epsilon$ for all $\epsilon > 0$ and consequently we have
\[
\mathbb{P}\{A_\epsilon\mid B_\epsilon\} = \frac{\mathbb{P}\{A_\epsilon \cap B_\epsilon\}}{\mathbb{P}\{B_\epsilon\}} = 1 \quad \forall~\epsilon > 0 .
\]
This shows that the regularity assumptions required in \cite{aizenman1993,aizenman2001} are violated. 
In the case $a < 0$ one can proceed analogously, but with a different choice of the sets $A_\epsilon$ and $B_\epsilon$.
%
\section{Proof of Theorem \ref{thm:apriori-trans}} \label{sec:proof}
Let $L \in \mathbb{N}$ be such that the cube $\Lambda^+ = [-L,L]^d \cap \mathbb{Z}^d$ contains
$\bigcup_{x \in \Lambda} \{k \in \mathbb{Z}^d \colon
u(x-k) \not = 0 \}$, which is the set of all lattice sites whose coupling
constant
influence the potential in $\Lambda$. 
Let $A:\ell^1 (\mathbb{Z}^d) \to \ell^1 (\mathbb{Z}^d)$ be the linear operator whose
coefficients
in the canonical orthonormal basis are $A(j,k) = u(j-k)$ for $j,k \in \mathbb{Z}^d$.
Since $u$ has compact support, the operator $A$ is bounded. Since
$\hat u$ does not vanish, $C_u = \lVert A^{-1} \rVert_1 < \infty$, see
\cite{veselic2008} for details. Moreover, there exists an invertible matrix
$A_{\Lambda^+} : \ell^1 (\Lambda^+) \to \ell^1 (\Lambda^+)$ satisfying
\begin{equation} \label{eq:circ}
 A_{\Lambda^+} (j,k) = u (j-k) \quad \text{for all $j\in \Lambda$ and $k \in
\Lambda^+$}
\end{equation}
and $\lVert A_{\Lambda^+}^{-1} \rVert \leq C_u$, see \cite[Proposition
5]{veselic2008}. We set $B_{\Lambda^+} = A_{\Lambda^+}^{-1}$.
From \cite{aizenman1993} and \cite{graf1994}, respectively, we infer 
that for $x,y \in \Lambda$ with $x \not = y$
\begin{equation} \label{eq:krein}
 \lvert G_{\omega,\Lambda} (z;x,x) \rvert = \frac{1}{\lvert V_\omega (x) -
\alpha\rvert}
\,\text{ and }\,
 \lvert G_{\omega,\Lambda} (z;x,y)\rvert \leq \frac{2}{\lvert V_\omega (x) -
\beta\rvert}  +
\frac{2}{\lvert V_\omega (y) - \gamma\rvert} .
\end{equation}
Here $\alpha,\beta \in \mathbb{C}$ are functions of
$V_\omega
(k)$, $k \in \Lambda \setminus \{x\}$, and $\gamma \in \mathbb{C}$ is a function of
$V_\omega (k)$, $k \in \Lambda \setminus \{y\}$. 
Set $\omega_{\Lambda^+} = (\omega_k)_{k \in
\Lambda^+}$, $k(\omega_{\Lambda^+}) = \prod_{k \in \Lambda^+} \rho
(\omega_k)$, $\ensuremath{\mathrm{d}} \omega_{\Lambda^+} = \prod_{k \in \Lambda^+} \ensuremath{\mathrm{d}} \omega_k$
and $\Omega_{\Lambda^+} = \times_{k \in \Lambda^+} \mathbb{R}$. Using the substitution
$\zeta_{\Lambda^+} = (\zeta_k)_{k \in \Lambda^+} = A_{\Lambda^+}
\omega_{\Lambda^+}$ and Eq. \eqref{eq:circ} we obtain
\begin{align*}
 S := \int_{\Omega_{\Lambda^+}} \frac{1}{|V_\omega (x) - \alpha|^s}
k(\omega_{\Lambda^+}) \ensuremath{\mathrm{d}} \omega_{\Lambda^+} 
= \int_{\Omega_{\Lambda^+}} \frac{1}{|\zeta_x - \alpha|^s} k(B_{\Lambda^+}
\zeta_{\Lambda^+}) \lvert \det B_{\Lambda^+} \rvert \ensuremath{\mathrm{d}} \zeta_{\Lambda^+} ,
\end{align*}
where $\ensuremath{\mathrm{d}} \zeta_{\Lambda^+} = \prod_{k \in \Lambda^+} \ensuremath{\mathrm{d}} \zeta_k$. Since by
construction $\zeta_k = V_\omega (k)$ for all $k \in \Lambda$, $\alpha$ is
now a function of $\zeta_k$, $k \in \Lambda \setminus \{x\}$. For non-negative
functions $g:\mathbb{R}\to\mathbb{R}$ with $g \in
W^{1,1} (\mathbb{R})$ and $\delta \in \mathbb{C}$ one has for all $\lambda > 0$ the estimate 
\begin{equation} \label{eq:graf}
 \int_\mathbb{R} \frac{g(\xi)}{\lvert \xi-\delta \rvert^s} \ensuremath{\mathrm{d}} \xi \leq \lambda^{-s}
\lVert g \rVert_{L^1} + \lVert g\rVert_\infty \frac{2\lambda^{1-s}}{1-s} 
\leq \lambda^{-s} \lVert g \rVert_{L^1} + \lVert g'\rVert_{L^1}
\frac{\lambda^{1-s}}{1-s} ,
\end{equation}
where the first inequality is due to \cite{graf1994} and the second is the
fundamental theorem of calculus. Set $\Omega_{\Lambda^+}^x = \times_{k \in
\Lambda^+ \setminus \{x\}} \mathbb{R}$. Using Fubini's theorem and Ineq. \eqref{eq:graf} we
obtain
\begin{align*}
 S &\leq \int\limits_{\Omega_{\Lambda^+}^x} \!\!\! \Bigl(
\lambda^{-s} \int\limits_\mathbb{R} k(B_{\Lambda^+} \zeta_{\Lambda^+}) \ensuremath{\mathrm{d}} \zeta_x +
\frac{\lambda^{1-s}}{1-s} \int\limits_\mathbb{R} \Bigl\lvert \frac{\partial}{\partial
\zeta_x}
k(B_{\Lambda^+} \zeta_{\Lambda^+}) \Bigr\rvert \ensuremath{\mathrm{d}} \zeta_x
\Bigr) |\det B_{\Lambda^+}| \!\!\!  \prod_{k \in \Lambda^+ \setminus \{x\}}
\!\!\!\!\! \ensuremath{\mathrm{d}}
\zeta_k \\
& = 
\lambda^{-s}  + \frac{\lambda^{1-s}}{1-s} \int_{\Omega_{\Lambda^+}} \Bigl\lvert
\frac{\partial}{\partial \zeta_x} k(B_{\Lambda^+} \zeta_{\Lambda^+}) \Bigr\rvert
 \lvert\det B_{\Lambda^+}\rvert  \ensuremath{\mathrm{d}} \zeta_{\Lambda^+} .
\end{align*}
We calculate the partial derivative by the product rule, substitute back into
original coordinates and obtain
\begin{align*}
S & \leq
\lambda^{-s}  + \frac{\lambda^{1-s}}{1-s} \int\limits_{\Omega_{\Lambda^+}}
\!\!\Bigl\lvert
\sum_{j \in \Lambda^+}  \rho'\bigl((B_{\Lambda^+} \zeta_{\Lambda^+})_j \bigr)
B_{\Lambda^+}(j,x) \!\! \prod_{\substack{k \in \Lambda^+ \\ k \not = j}} \!\!
\rho\bigl(
(B_{\Lambda^+} \zeta_{\Lambda^+})_k \bigr)  \Bigr\rvert
 \lvert\det B_{\Lambda^+}\rvert  \ensuremath{\mathrm{d}} \zeta_{\Lambda^+} \\
&\leq \lambda^{-s}  + \frac{\lambda^{1-s}}{1-s}  \sum_{j \in \Lambda^+}
\bigl\lvert B_{\Lambda^+}(j,x) \bigr\rvert \int_{\Omega_{\Lambda^+}}
 \bigl\lvert \rho' (\omega_j ) \bigr\rvert
\prod_{\substack{k \in \Lambda^+ \\ k \not = j}} \rho ( \omega_k )  
  \ensuremath{\mathrm{d}} \omega_{\Lambda^+} \\
&= \lambda^{-s}  + \frac{\lambda^{1-s}}{1-s} \lVert B_{\Lambda^+} \rVert_1
\lVert \rho' \rVert_{L^1} \leq \lambda^{-s}  + \frac{\lambda^{1-s}}{1-s} C_u
\lVert \rho' \rVert_{L^1}.
\end{align*}
If we choose $\lambda = s/(C_u \lVert \rho' \rVert_{L^1})$ we obtain $S \leq
C_u^s \lVert \rho' \rVert_{L^1}^s s^{-s}/(1-s)$. Thus we have shown $\mathbb{E}
\bigl\{| G_{\omega,\Lambda} (z;x,x) |^s \bigr\} \leq C_u^s
\Vert\rho'\Vert_{L^1}^s \frac{s^{-s}}{1-s}$. Analogously we obtain for $x \not =
y$ the estimate $\mathbb{E} \bigl\{ | G_{\omega , \Lambda} (z;x,y) |^s \bigr\} \leq
2^{s+1}C_u^s \Vert\rho'\Vert_{L^1}^s \frac{s^{-s}}{1-s}$. The additional factor
$2^{s+1}$ arises in the case $x \not = y$, since in Eq. \eqref{eq:krein} there
are two summands of the type $S$, and each summand has an additional factor $2$
in its numerator. \qed

\end{document}